\documentclass[11pt, reqno, english]{amsart}  
\usepackage[utf8]{inputenc}
\usepackage[T1]{fontenc}
\usepackage{amsmath,amsthm}
\usepackage{amsfonts,amssymb}
\usepackage{url}
\usepackage{mathtools}  
\usepackage[colorlinks=true,urlcolor=blue,linkcolor=red,citecolor=magenta]{hyperref}
\usepackage{enumerate, paralist}
\usepackage{tikz}

\usepackage[left=1in,right=1in,top=1in,bottom=1in]{geometry}
\numberwithin{equation}{section}

\theoremstyle{plain}
\newtheorem{theorem}{Theorem}[section]

\newtheorem{proposition}[theorem]{Proposition}

\theoremstyle{definition}

\newcommand{\R}{\mathbb{R}}

\newcommand{\Z}{\mathbb{Z}}

\newcommand{\F}{\mathcal{F}}

\DeclareMathOperator{\Conv}{\mathrm{Conv}}

\DeclareMathOperator{\Relint}{\mathrm{Relint}}

\title{A Non-face Characterization of Spheres on Few Vertices}

\author[Huang]{Shuai Huang}
\address[SH]{Dept.\ Math.\, Bard College, Annandale-on-Hudson, NY 12504, USA}
\email{eh6041@bard.edu} 

\author[Miller]{Jasper Miller}
\address[JM]{Dept.\ Math.\, Bard College, Annandale-on-Hudson, NY 12504, USA}
\email{jm5305@bard.edu} 

\author[Rose]{Daniel Rose-Levine}
\address[DR]{Dept.\ Math.\, Bard College, Annandale-on-Hudson, NY 12504, USA}
\email{dr6048@bard.edu} 

\author[Simon]{Steven Simon}
\address[SS]{Dept.\ Math.\, Bard College, Annandale-on-Hudson, NY 12504, USA}
\email{ssimon@bard.edu} 

\begin{document}

\begin{abstract} We prove a relatively simple combinatorial characterization of simplicial $d$-spheres on $d+4$ vertices. Our criteria are given in terms of the intersection patterns of a simplicial complex's family of minimal non-faces. Namely, let $\Sigma$ be a simplicial complex on $d+4$ vertices and let $\F$ be its family of minimal non-faces. Then $\Sigma$ is a $d$-sphere if and only if $|\F|=n\geq 3$ is odd and there is an ordering $A_0,\ldots, A_{n-1}$ of the minimal non-faces, indices taken modulo $n$, such that successive $A_i$ are disjoint and the alternating $\frac{(n-1)}{2}$-fold intersections $A_i\cap A_{i+2} \cap A_{i+4} \cap \cdots \cap A_{i+n-3}$ partition the vertex set.   
\end{abstract}

\maketitle

\section{Introduction}

A basic problem in discrete geometry, for which very few general results are known, seeks necessary and sufficient combinatorial criteria to determine when an arbitrary simplicial complex $\Sigma$ is a simplicial $d$-sphere. Computationally, deciding when $\Sigma$ is a simplicial $d$-sphere is similarly difficult: while possible when $d\le3$~\cite{Ru95}, it is impossible for all $d\geq 5$~\cite{No55, CL06} and remains open when $d=4$~\cite{Po14}. Here we show that if one restricts to simplicial complexes on relatively few vertices (i.e., at most $d+4$), however, then surprisingly simple combinatorial characterizations of simplicial $d$-spheres exist. We note that all simplicial $d$-spheres on at most $d+4$ vertices are polytopal by a theorem of Mani~\cite{Ma72}. 

Our characterization of spheres on few vertices is given in terms of the simplicial complex's family $\F$ of minimal non-faces. Recall that for a simplicial complex $\Sigma$ with vertex set $[m]=\{1,\ldots,m\}$, its \emph{non-faces} are the elements of $2^{[m]}$ which are not elements of the complex, and that a non-face is \emph{minimal} if any of its proper subsets is a face of $\Sigma$. Note that any $A\in \F$ has size at least two because every singleton of $[m]$ is a face of $\Sigma$ by assumption. It is easily seen that $\F$ completely determines $\Sigma$; see the discussion immediately preceding the examples of Theorem~\ref{thm:d+4} below. We remark that $\F$ is the set of minimal generators of the Stanley-Reisner ideal of $\Sigma$, and also that the intersection patterns of the minimal non-faces serve as the criteria of powerful non-embeddability theorems (see, e.g.,~\cite{Ma08, Sa90, Sa91}). 

Trivially, a simplicial complex $\Sigma$ with vertex set $[d+2]$ is a $d$-sphere if and only if it is the boundary complex of the $(d+1)$-dimensional simplex, or equivalently that $\F=\{[d+2]\}$ is the trivial partition of $[d+2]$. For simplicial complexes $\Sigma$ with vertex set $[d+3]$, it is implicit from~\cite{FS23} (see Theorem~1 and the proof of Lemma 8) that $\Sigma$ is a $d$-sphere if and only if $\F$ is a partition of $[d+3]$ by two sets of order at least two.

Now suppose that $\Sigma$ is a simplicial complex with vertex set $[d+4]$ and let $n=|\F|$. We say that $\F$ is an $n$-\emph{cycle} if there is an ordering $A_0,\ldots, A_{n-1}$ of the elements of $\F$, indices taken in $\Z_n$, such that successive $A_i$ are disjoint, i.e, $A_i\cap A_{i+1}=\emptyset$ for all $i\in \Z_n$. We call such an ordering a \emph{cyclic ordering} of the elements of $\F$. Now suppose that $n\geq 3$ is odd and that $\F$ is an $n$-cycle with some cyclic ordering $A_0,\ldots, A_{n-1}$. Then for each $i\in \Z_n$ one may consider the alternating $(\frac{n-1}{2})$-fold intersection \[B_i:=A_i\cap A_{i+2}\cap A_{i+4}\cap\cdots\cap A_{i+n-3}.\] Thus the $B_i$ represent the intersections $A_{j_1}\cap \cdots \cap A_{j_k}$ of maximum length for pairwise distinct $A_{j_i}$ which are potentially non-empty. We say that $\F$ is a \emph{maximum odd cycle} if $n\geq 3$ is odd and if $\F$ is an $n$-cycle for which there is some cyclic ordering $A_0,\ldots, A_{n-1}$ of the elements of $\F$ such that $B_0,\ldots, B_{n-1}$ partition $[d+4]$ (thus $n\leq d+4$). One then has the following characterization of simplicial $d$-spheres on $d+4$ vertices:

 \begin{theorem}
 \label{thm:d+4} Let $d\geq 1$ be an integer, let $\Sigma$ be a simplicial complex with vertex set $[d+4]$, and let $\F$ be its family of minimal non-faces. Then $\Sigma$ is a simplicial $d$-sphere if and only if $\F$ is a maximum odd cycle. 
\end{theorem}

To our knowledge Theorem~\ref{thm:d+4} does not follow from either Mani's theorem~\cite{Ma72} above or from the extensive known results concerning $d$-polytopes with $d+3$ vertices (see, e.g.,~\cite{Gr03, Mc74, Mc78, DDS04}). Indeed the sufficiency of the criteria of Theorem~\ref{thm:d+4} is independent of these. The necessity of these criteria does rely on~\cite{Ma72}, however, as well as basic Gale duality applied to such polytopes.

Before giving examples of Theorem~\ref{thm:d+4} we first set some notation. Namely, to  any family $\F$ of subsets of $[m]$, each of order at least two, we assign the simplicial complex $\Sigma(\F)$ defined by \[\Sigma(\F)=\{A\subseteq [m] \mid B\notin \F\,\, \text{if}\,\, B\subseteq A\}.\] Thus $\Sigma(\F)$ is the simplicial complex defined by setting $A\in \Sigma(\F)$ if and only if $A$ does not contain any element of $\F$ as a subset. It is easily observed that the vertex set of $\Sigma(\F)$ is $[m]$ and that the set of minimal non-faces of $\Sigma(\F)$ is precisely $\F$. Moreover, for any simplicial complex $\Sigma$ one has $\Sigma=\Sigma(\F)$ where $\F$ is the family of minimal non-faces of $\Sigma$.

As a first example of Theorem~\ref{thm:d+4}, consider when $n=3$. For any subset $A\subseteq [d+4]$ with $|A|\geq 2$, let $\Delta_A$ denote the $(|A|-1)$-dimensional simplex with vertex set $A$ and let $\partial \Delta_A$ denote its boundary complex. Thus $\partial \Delta_A$ is an $(|A|-2)$-dimensional sphere.  Now let $\F=\{A_0,A_1,A_2\}$ be any partition of $[d+4]$ with $|A_i|\geq 2$ for all $i$. As $A_i=B_i$ for all $i$, $\F$ is a maximum $3$-cycle for $\Sigma:=\Sigma(\F)$. One has that $\Sigma=\partial \Delta_{A_1}\ast \partial \Delta_{A_2}\ast\partial \Delta_{A_3}$ is the join of the $\partial \Delta_{A_i}$ and so $\Sigma$ is a $d$-sphere. In particular, when $d=2$ the family $\F=\{\{1,2\},\{3,4\},\{5,6\}\}$ determines an octahedron. As a second example, note that when $d=1$ the only possibility for a maximum odd cycle is when $n=5$. Letting  \[\F=\{\{1,4\},\{2,5\},\{1,3\},\{2,4\},\{3,5\}\},\] consider $\Sigma:=\Sigma(\F)$. As $B_i=A_i\cap A_{i+2}=\{i+1\}$ for all $i\in \{0,1,2,3,4\}$, $\F$ is a maximum $5$-cycle for $\Sigma$. It is readily verified that $\Sigma$ is the $5$-cycle graph with edges $\{1,2\},\{2,3\},\{3,4\},\{4,5\},\{1,5\}$, which is a circle.

\section{Preliminaries on Gale Duality}
\label{sec:prelim}

Our proof of Theorem~\ref{thm:d+4} is based on Gale duality~\cite{Ga56}, an essential tool in polyhedral combinatorics. We briefly recall its construction and the fundamental properties we need for our purposes; the reader is referred to ~\cite {Gr03,Ma02} for a more detailed discussion.

 Let $x_1,\ldots, x_n$ be a sequence of affinely spanning points in $\R^d$ ($n\geq d+1$). Correspondingly, the $(d+1)\times n$ matrix $A=\begin{bmatrix} x_1\ldots x_n\\
                     1 \ldots 1\end{bmatrix}$
whose columns are the $(x_i,1)\in \R^{d+1}$ has full rank. Thus $\dim (\ker A) =n-d-1$. Fixing a basis $b_1,\ldots, b_{n-d-1}$ for $\ker A$, one defines the \emph{Gale transform} of $x_1,\ldots, x_n$ to be the row vectors $y_1,\ldots, y_n\in \R^{n-d-1}$ of the $n\times (n-d-1)$ matrix $B$ whose columns are the $b_i$. It is easily seen that (1) the $y_i$ linearly span $\R^{n-d-1}$ and (2) $\sum_{i=1}^n y_i=0$. Conversely, if $y_1,\ldots, y_n$ are a sequence of vectors in $\R^{n-d-1}$ satisfying (1) and (2), then it follows that there exist affinely spanning points $x_1,\ldots, x_n$ in $\R^d$ whose Gale transform are the $y_i$. 

A central property of the Gale transform is that it links affinely dependences among the $x_i$ in $\R^d$ with separations of the $y_i$ by linear hyperplanes in $\R^{n-d-1}$. Let $\langle u, w\rangle=\sum_{i=1}^n u_iw_i$ denote the usual inner product of two vectors $u$ and $w$ in $\R^n$.

\begin{proposition}
\label{prop: Gale dependence} Let $y_1,\ldots, y_n \in \R^{n-d-1}$ be the Gale transform of a sequence $x_1,\ldots, x_n$ of affinely spanning points in $\R^d$. If $\lambda=(\lambda_1,\ldots, \lambda_n)$ is a non-zero vector in $\R^n$ such that (1) $\sum_{i=1}^n\lambda_ix_i=0$ and (2) $\sum_{i=1}^n\lambda_i=0$, then there exists a non-zero vector $\alpha$ in $\R^{n-d-1}$ such that $\langle \alpha, y_i \rangle =\lambda_i$ for all $i\in [n]$. Conversely, given a non-zero vector $\alpha$ in $\R^{n-d-1}$, define $\lambda_i=\langle \alpha, y_i \rangle$ for all $i\in [n]$. Then $\lambda=(\lambda_1,\ldots, \lambda_n)\in \R^n$ is non-zero, $\sum_{i=1}^n\lambda_ix_i=0$, and $\sum_{i=1}^n\lambda_i=0$. 
\end{proposition}

Given $x_i$ as above, their \emph{Gale diagram} $\widehat{y}_1,\ldots, \widehat{y}_n$ is defined by setting $\widehat{y_i}=0$ if $y_i=0$ and $\widehat{y}_i=y_i/\|y_i\|$ otherwise. Thus the Gale diagram lies in $S^{n-d-2}\cup \{0\}$. The central feature of the Gale diagram is that it completely determines the face lattice of a $d$-polytope $P:=\Conv(\{x_1,\ldots, x_n\})$ (see Theorem~1 on pg. 88 and comments on pg. 89 of~\cite{Gr03}). As usual, we let $\Relint(S)$ denotes the relative interior of a subset $S$ of a convex set.

\begin{theorem}
\label{thm:Gale face}
Let $P=\Conv(\{x_1,\ldots, x_n\})$ be a $d$-polytope in $\R^d$ with vertices $x_1,\ldots, x_n$ and let $\widehat{y}_1,\ldots, \widehat{y}_n\in S^{n-d-2}\cup \{0\}$ be the corresponding Gale diagram of the $x_i$. Suppose that  $A\subset [n]$. Then $\Conv(\{x_i\mid i \in A\})$ is a proper face of $P$ if and only if $0\in \Relint(\Conv (\{\widehat{y}_i\mid i \notin A\}))$.
\end{theorem}

\section{Proof of Theorem~\ref{thm:d+4}}

We now prove Theorem~\ref{thm:d+4}.

\subsection{Proof of Sufficiency} We first prove the sufficiency of the criteria in Theorem~\ref{thm:d+4}.

\begin{proof} Let $k\geq 1$ and let $\Sigma$ be a simplicial complex with vertex set $[d+4]$ whose family $\F=\{A_0,\ldots, A_{2k}\}$ of minimal non-faces is a maximum $(2k+1)$-cycle, indices taken in $\Z_{2k+1}$. We have that $\Sigma=\Sigma(\F)$; we show that $\Sigma(\F)$ can be identified with the boundary complex $\partial P$ of a simplicial $(d+1)$-polytope $P$ in $\R^{d+1}$ with $d+4$ vertices. Then $\Sigma$ is a $d$-sphere.

For each $i\in \Z_{2k+1}$ we have $B_i=A_i\cap A_{i+2} \cap A_{i+4}\cap \cdots \cap A_{i+2k-2}$. Thus the $B_i$ are obtained by cyclically shifting the indices of $B_0=A_0\cap A_2\cap A_4\cap \cdots \cap A_{2k-2}$, from which it follows easily that $B_0\cup B_{-2}\cup B_{-4}\cup \cdots \cup B_{-2(k-1)}\subseteq A_0$. On the other hand, $A_0\cap B_i=\emptyset$ if $i\notin \{0,2,\ldots, 2k-2\}$ because then $0$ is adjacent to one of $i,i+2\ldots, i+2k-2$. Since the $B_i$ partition $[d+4]$ we conclude that $A_0=B_0\cup B_{-2}\cup B_{-4}\cup \cdots \cup B_{-2(k-1)}$, and an identical argument shows that 
\begin{equation}\label{eqn: A_i} A_i=B_{i}\cup B_{i-2}\cup B_{i-4}\cup \cdots \cup B_{i-2(k-1)}\end{equation} for all $i\in \Z_{2k+1}$. 

In order to apply Gale duality, we view the $B_i$ as equally spaced points on the unit circle $S^1\subset \R^2$, counted with multiplicity according to the size of $B_i$, in the cyclic order \begin{equation}\label{eqn: B_i} B_0, B_{-2}, B_{-4},\ldots, B_{-4k}. \end{equation} Noting that $B_{-4k-2}=B_0$, we see from equation~(\ref{eqn: A_i}) that each element of $\F$ is represented by a sequence of $k$ consecutive $B_j$'s in the order given by~(\ref{eqn: B_i}), and conversely; see Figure~1 below for an example when $k=2$. Explicitly, let $v_0,\ldots, v_{2k}$ be the vertices in cyclic order of a regular $(2k+1)$-gon which is inscribed in $S^1$. For each $j\in \Z_{2k+1}$, let $y_i=\frac{1}{|B_{-2j}|}v_j\in \R^2$ provided $i\in B_{-2j}$.  We have that $y_1,\ldots, y_{d+4}$ linearly span $\R^2$, and since $\sum_{i\in [d+4]}y_i=\sum_{j=0}^{2k-1}v_j=0$ they are the Gale transform of some sequence $x_1,\ldots, x_{d+4}$ in $\R^{d+1}$. Now let $\widehat{y}_1,\ldots, \widehat{y}_{d+4}\in S^1$ be the Gale diagram of the $x_i$. Thus for each $j\in \Z_{2k+1}$ one has $\widehat{y}_i=v_j$ if and only if $i\in B_{-2j}$, and so each $v_j$ corresponds to $B_{-2j}$. 

\begin{figure}[ht]
\centering
\label{fig:fig}
 \noindent
\begin{minipage}[h]{0.65\textwidth}\hspace*{.65in}
  \begin{tikzpicture}[scale=2.2]

    \node[fill=black, circle, inner sep=1.5pt, label={[xshift=4pt, yshift=1pt]above:$B_3$}] at (72:1) {};
    \node[fill=black, circle, inner sep=1.5pt, label=above left:$B_1$] at (144:1) {};
    \node[fill=black, circle, inner sep=1.5pt, label=below left:$B_4$] at (216:1) {};
    \node[fill=black, circle, inner sep=1.5pt, label={[xshift=4pt, yshift=-2pt]below:$B_2$}] at (288:1) {};
    \node[fill=black, circle, inner sep=1.5pt, label=right:$B_0$] at (0:1) {};
    \node[fill=black, circle, inner sep=1.5pt, label={[xshift=4pt]below:$O$}] at (0,0) {};

    \draw[bend right=30] (72:1) to (144:1);
    \draw[bend right=30] (144:1) to (216:1);
    \draw[bend right=30] (216:1) to (288:1);
    \draw[bend right=30] (288:1) to (0:1);
    \draw[bend right=30] (0:1) to (72:1);

    \def\i{6}
    \def\j{60}
    \coordinate (P1) at ({72+\i}:1.1);
    \coordinate (P2) at ({144+\i}:1.1);
    \coordinate (P3) at ({216+\i}:1.1);
    \coordinate (P4) at ({288+\i}:1.1);
    \coordinate (P5) at ({0+\i}:1.1);
    \coordinate (P1_) at ({72-\i}:1.1);
    \coordinate (P2_) at ({144-\i}:1.1);
    \coordinate (P3_) at ({216-\i}:1.1);
    \coordinate (P4_) at ({288-\i}:1.1);
    \coordinate (P5_) at ({0-\i}:1.1);

    \draw[blue, bend right=\j] (P1) to node[midway, above, blue] {$A_3$} (P2_);
    \draw[blue, bend right=\j] (P2) to node[midway, left, blue] {$A_1$} (P3_);
    \draw[blue, bend right=\j] (P3) to node[midway, below, blue] {$A_4$} (P4_);
    \draw[blue, bend right=\j] (P4) to node[midway, right=2pt, blue] {$A_2$} (P5_);
    \draw[blue, bend right=\j] (P5) to node[midway, right=2pt, blue] {$A_0$} (P1_);
    
  \end{tikzpicture}
\end{minipage}

\caption{The case $k=2$}
\end{figure}
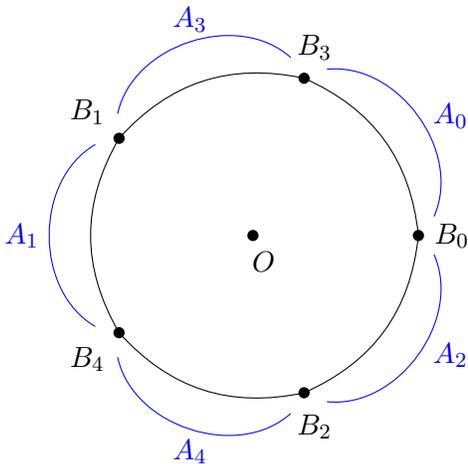

Now let $P=\Conv(\{x_1,\ldots, x_{d+4}\})$, which is a $(d+1)$-polytope. We first show that the $x_i$ are the vertices of $P$, i.e., that the $x_i$ are extremal. This follows from Proposition~\ref{prop: Gale dependence}. For if $x_1$, say, were in the convex hull of the remaining $x_j$ then there would exist a line $L$ in $\R^2$ through the origin which contains only $y_1$ in an open half-space $H^+$ determined by $L$. This is the case if and only if $\widehat{y_1}$ lies in $H^+$. We now  consider the cases $k=1$ and $k\geq 2$ separately. In the first case we have that $B_i=A_i$ for all $i$, and we must have $|A_i|\geq 2$; this is because the emptyset is a face of $\Sigma$, as is any singleton. We have thus reached a contradiction since $\widehat{y_1}$ occurs with multiplicity at least two. For $k\geq 2$ we also reach a contradiction since $H^+$ must contain at least $k$ of the $v_j$'s and hence of the $y_j$'s. Thus the $x_i$ are extremal.

Next we show that any $A\subset [d+4]$ is a face of $\Sigma(\F)$ if and only if $\Conv(\{x_i\mid i \in A\})$ is a proper face of $P$. Thus $P$ is simplicial and $\Sigma(\F)$ can be identified with the boundary complex $\partial P$. So let $A\in \Sigma(\F)$. To show that $\Conv(\{x_i\mid i \in A\})$ is a proper face of $P$, by Theorem~\ref{thm:Gale face} we show that $0\in \Relint(\Conv(\{\widehat{y_i}\mid i\notin A\}))$. Now the latter holds if and only if the $\widehat{y}_i$ with $i\notin A$ do not all arise from some sequence of $k+1$ consecutive $B_i$'s in the order prescribed by~(\ref{eqn: B_i}), or more formally if for any choice $v_i,\ldots, v_{i+k+1}$ of $k+1$ consecutive vertices there is some $\widehat{y}_{j_0}$ with $j_0\notin A$ such that $\widehat{y}_{j_0}\notin V:=\{v_i,\ldots, v_{i+k+1}\}$. So let $v_i,\ldots, v_{i+k+1}$ be given and consider their complementary vertices on $S^1$. These $k$ consecutive vertices correspond to the sequence $B_j, B_{j-2},\ldots, B_{j-2(k-1)}$ representing some $A_j\in \F$. Since $A\in \Sigma(\F)$ we have $A_j\not\subseteq A$. Thus there exists some $j_0\in A_j$ such that $j_0\notin A$. As $j_0\in A_j$, $j_0$ lies in some $B_{j-2\ell}$ for some $\ell\in\{0,\ldots, k-1\}$. We therefore have that $\widehat{y}_{j_0}\notin V$ and thus that $\Conv(\{x_i\mid i \in A\})$ is a face of $P$. Conversely, let $A\subset [d+4]$ and suppose that $\Conv(\{x_i\mid i \in A\})$ is a proper face of $P$. Let $A_i\in \F$ be given, and consider the corresponding arc of $k$ consecutive vertices given by the $B_j$'s which compose $A_i$. Let $V$ be the set of complementary vertices of this arc. Since $0\in \Relint(\Conv(\{\widehat{y}_i\mid i \notin A\}))$, there is some $i_0\notin A$ such that $\widehat{y}_{i_0}$ does not lie in $V$. Thus $i_0$ lies in one of the $B_j$'s composing $A_i$, so $i_0\in A_i$. Thus $A_i\not\subseteq A$ for any $i\in \Z_{2k+1}$ and so $A\in \Sigma(\F)$.\end{proof}

\subsection{Proof of Necessity} We now prove the necessity of the criteria of Theorem~\ref{thm:d+4}. 

\begin{proof} Suppose that $\Sigma$ is a simplicial sphere with vertex set $[d+4]$. By Mani's theorem~\cite{Ma72}, $\Sigma$ can be identified with the boundary complex of a $(d+1)$-dimensional simplicial polytope $P$ in $\R^{d+1}$ on $d+4$ vertices, say with vertices $x_1,\ldots, x_{d+4}$. Let $\widehat{y}_1,\ldots, \widehat{y}_{d+4}$ be their Gale diagram in $S^1\cup \{0\}$. Since $P$ is simplicial, we have by Theorem~\ref{thm:Gale face} that $\widehat{y}_j\neq 0$ for all $j$ and that no two $\widehat{y}_j$ lie on antipodal points of $S^1$. It follows (see ~\cite[Pg. 111]{Gr03}) that the Gale diagram can be represented by $2k+1$ equally spaced points $v_0,\ldots, v_{2k}$ on the unit circle, counted with multiplicity, for some $k\geq 1$. One now only need reverse the steps above. Namely, for all $j\in \Z_{2k+1}$ define $B_{-2j}=\{i\in [d+4]\mid \widehat{y}_i=v_j\}$. Thus $B_0,\ldots, B_{2k}$ partition $[d+4]$. Next, define $A_0,\ldots, A_{2k}$ by equation~(\ref{eqn: A_i}) and let $\F=\{A_0,\ldots, A_{2k}\}$. The identical argument given above shows that $A\subset [d+4]$ is a face of $\Sigma(\F)$ if and only if $\Conv(\{x_i\mid i \in A\})$ is a proper face of $P$. As this holds for $\Sigma$ as well we conclude that $\Sigma=\Sigma(\F)$ and therefore that $\F$ is the family of minimal non-faces of $\Sigma$. 

Finally, we show that $\F$ is a maximum $(2k+1)$-cycle. However, it follows easily from the definition of the $A_i$ and that the $B_j$ are pairwise disjoint that $A_i\cap A_{i+1}=\emptyset$ for all $i\in \Z_{2k+1}$. Thus $\F$ is a $(2k+1)$-cycle. Similarly, one sees for each $i\in \Z_{2k+1}$ that $A_i \cap A_{i+2}=B_i\cap B_{i-2}\cap \cdots \cap B_{i-2(k-2)}$, and continuing that  $A_i\cap A_{i+2}\cap \cdots \cap A_{i+2(k-1)}=B_i$. As the $B_j$ partition $[d+4]$ this completes the proof. \end{proof}

\section*{Acknowledgment}

The authors are grateful for the support provided by the 2024 Bard Science Research Institute. They also thank Florian Frick for helpful comments on the manuscript.


\begin{thebibliography}{}


\bibitem{CL06} A. V. Chernavsky and V.P. Leksine. Unrecognizability of manifolds. \emph{Ann. Pure Appl. Logic} Vol.  141 No. 3 (2006) 325--335.

\bibitem{DDS04} M. Deza, M. Dutour, and M. Shtogrin. On simplicial and cubical complexes with short links. \emph{Israel J. Math}, Vol. 144 No. 1 (2004) 109--124.

\bibitem{FS23} F. Frick, M. Hu, V. Scheel, and S. Simon. Embedding Dimensions of Simplicial Complexes on Few Vertices. \emph{Ann. Comb.} Vol. 27 (2023) 993--1003.

\bibitem{Ga56} D. Gale. Neighboring vertices on a convex polyhedron, \textit{Linear inequalities and related system}, Annals of Mathematics Studies, No. 38, Princeton University Press, Princeton, N.J. (1956) 255--263.

\bibitem{Gr03} B. Gr\"unbaum. \emph{Convex Polytopes}, Second Edition, Graduate Texts in Mathematics Vol. 221, Springer--Verlag, New York (2003). 

\bibitem{Ma72} P. Mani. Spheres with few vertices. \emph{J. Comb. Theory, Ser. A}, Vol. 13 No.3 (1972) 346--352.

\bibitem{Ma02} J. Matou\v sek. \textit{Lectures on Discrete Geometry}, Graduate Texts in Mathematics Vol. 212, Springer--Verlag, New York (2002).

\bibitem{Ma08} J. Matou\v sek. \emph{Using the Borsuk--Ulam theorem: Lectures on topological methods in combinatorics and geometry.} Springer-Verlag (2008).

\bibitem{Mc74} P. McMullen. The number of neighbourly d-polytopes with $d+3$ vertices. \emph{Mathematica}, Vol. 21 No. 1 (1974) 26--31.

\bibitem{Mc78} P. McMullen. Transforms, diagrams and representations. \emph{Contributions to Geometry}, J. Tolke and J. M. Wills, eds., Birkhauser, Basel (1978) 92--130.

\bibitem{No55} P.S. Novikov. On the algorithmic insolvability of the word problem in group theory. \emph{Izdat. Akad. Nauk SSSR, Moscow: Trudy Mat}, p. 44. Inst. Steklov. no (1955).

\bibitem{Po14} B. Poonen. Undecidable problems: a sampler. \emph{Interpreting G\"odel: Critical Essays}, J. Kennedy, ed., Cambridge University Press (2014) 211-241.

\bibitem{Ru95} J.H. Rubinstein. An algorithm to recognize the 3-sphere. \emph{Proceedings of the International Congress of Mathematicians}, Vol. 1 and 2, Birkh\"auser, Basel (1995) 60--611. 

\bibitem{Sa90} K. S. Sarkaria. A generalized Kneser conjecture, \textit{J. Combin. Theory, Ser. B}, Vol. 49 No. 2 (1990) 236--240.

\bibitem{Sa91} K. S. Sarkaria, A generalized van Kampen--Flores theorem, \textit{Proc. Amer. Math. Soc.}, Vol. 11 (1991) 559--565.



\end{thebibliography}
\end{document}